\documentclass{article}
\catcode`\@=11
\newtheorem{theorem}{Theorem}
\newtheorem{lemma}{Lemma}
\newtheorem{proposition}{Proposition}

\newtheorem{example}{Example}
\newtheorem{remark}{Remark}

\newtheorem{corollary}{Corollary}
\def\demo{\noindent{\bf Proof .-}}
\def\section{\@startsection {section}{1}{\z@}{-3.5ex plus -1ex
minus-.2ex}{2.3ex plus .2ex}{\normalsize\bf}}
\def\codim{{\rm codim}\,}

\def\bz{\hbox{\it Z\hskip -4pt Z}}

\def\bc{\hbox{\,\it l\hskip -5.5pt C\/}}

\newcommand{\het}{H_{\rm et}}
\newcommand{\hc}{H_{\rm c}}

\begin{document}
\begin{center}
{\Large\bf \textsc{On the arithmetical rank of a special class of minimal varieties}}
\end{center}
\vskip.5truecm
\begin{center}
{Margherita Barile\footnote{Partially supported by the Italian Ministry of Education, University and Research.}\\ Dipartimento di Matematica, Universit\`{a} di Bari,Via E. Orabona 4,\\70125 Bari, Italy\\barile@dm.uniba.it}
\end{center}
\vskip1truecm
\noindent
{\bf Abstract} We study the arithmetical ranks and the cohomological dimensions of an infinite class of Cohen-Macaulay varieties of minimal degree. Among these we find, on the one hand, infinitely many set-theoretic complete intersections, on the other hand examples where the arithmetical rank is arbitrarily greater than the codimension. 

\section*{Introduction}
Let $K$ be an algebraically closed field, and let $R$ be a polynomial ring  in $N$ indeterminates over $K$. Let $I$ be a proper reduced ideal of $R$ and consider the variety $V(I)$ defined in the affine space $K^N$ (or in the projective space ${\bf P}_K^{N-1}$, if $I$ is homogeneous and different from the maximal irrelevant ideal) by the vanishing of all polynomials in $I$. By Hilbert's Basissatz there are  finitely many polynomials $F_1,\dots, F_r\in R$ such that $V(I)$ is defined by the equations $F_1=\cdots= F_r=0$. By Hilbert's Nullstellensatz this is equivalent to the ideal-theoretic condition
$$I=\sqrt{(F_1,\dots, F_r)}.$$
\noindent
Suppose $r$ is minimal with respect to this property. It is well known that $\codim V(I)\leq r$. If equality holds, $V(I)$ is called a {\it set-theoretic complete intersection} on $F_1,\dots, F_r$. \newline
Exhibiting significant examples of set-theoretic complete intersections (or, more generally, determining the minimum number of equations defining given varieties, the so-called {\it arithmetical rank}, denoted ara, of their defining ideals) is one of the hardest problems in algebraic geometry. In \cite{Ba2} we already determined infinitely many set-theoretic complete intersections among the Cohen-Macaulay varieties of minimal degree which were classified geometrically by Bertini \cite{B}, Del Pezzo \cite{DP}, Harris \cite{Ha} and Xamb\'o \cite{X}, and whose defining ideals were determined in an explicit algorithmic way in \cite{BM2}. In this paper we present a new class of minimal varieties, where the gap between the arithmetical rank and the codimension can be arbitrarily high. It includes an infinite set of set-theoretic complete intersections. For the arithmetical ranks of the complementary set of varieties we determine a lower bound (given by \'etale cohomology) and an upper bound (resulting from the computation of an explicit set of defining equations) that  only differ by one: the equality between the lower bound and the actual value of the arithmetical rank is shown in few  special cases. We also determine the cohomological dimensions of the defining ideals of each of these varieties. This invariant, in general, also provides a lower bound for the arithmetical rank, and the cases where it is known to be smaller are rare. Those which were found so far are the determinantal and Pfaffian ideals considered in \cite{BS} and in \cite{Ba1}: there the strict inequality holds in all positive characteristics. We prove that the same is true for the minimal varieties investigated in the present paper that are not set-theoretic complete intersections.\newline
Some crucial results on arithmetical ranks and cohomological dimensions are due to Bruns et al. and are quoted from \cite{BS} and \cite{BV}.
\section{Preliminaries}
For all integers $s\geq2$ and $t\geq1$ consider the two-row matrix
$$A_{s,t}=\left(
\begin{array}{ccccccccccccc}
x_1&x_2&\cdots&x_s&\vline\ \vline&y_0&\vline\ \vline&y_1&\vline\ \vline&\cdots&\vline\ \vline&y_{t-1}\\
x_{s+1}&x_{s+2}&\cdots&x_{2s}&\vline\ \vline&z_1&\vline\ \vline&z_2&\vline\ \vline&\cdots&\vline\ \vline&z_t\\
\end{array}\right),
$$
\noindent
where $x_1, x_2,\dots, x_{2s},y_0, y_1, \dots, y_{t-1}, z_1, z_2,\dots, z_t$ are $N$ indeterminates over $K$. We assume that they are pairwise distinct, possibly with the following exception: we can have $x_{2s}=y_0$ or $z_i=y_j$ for some indices $i$ and $j$ such that $1\leq i\leq j\leq t-1$, but no entry appears more than twice in $A_{s,t}$. We have the least possible number of indeterminates if $x_{2s}=y_0$ and $z_i=y_i$ for  $1=1,\dots, t-1$, in which case $N=2s+t$, and the matrix takes the following form:
$$\bar A_{s,t}=\left(
\begin{array}{ccccccccccccc}
x_1&x_2&\cdots&x_s&\vline\ \vline&x_{2s}&\vline\ \vline&y_1&\vline\ \vline&\cdots&\vline\ \vline&y_{t-1}\\
x_{s+1}&x_{s+2}&\cdots&x_{2s}&\vline\ \vline&y_1&\vline\ \vline&y_2&\vline\ \vline&\cdots&\vline\ \vline&y_t\\
\end{array}\right).
$$
\noindent
If the indeterminates are pairwise distinct, then $N=2s+2t$.  The matrix $A_{s,t}$ belongs to the class of so-called {\it barred matrices} introduced in \cite{BM1} and can be associated with the ideal $J_{s,t}$ of $R=K[x_1,x_2,\dots, x_{2s},y_0, y_1,\dots, y_{t-1}, z_1, z_2,\dots, z_t]$ generated by the union of 
\begin{list}{}{}
\item{(I)} the set ${\cal M}$ of two-minors of the submatrix of $A_{s,t}$ formed by the first $s$ columns (the so-called first {\it  big block});
\item{(II)} the set of products $x_iz_j$, with $1\leq i\leq s$ and $1\leq j\leq t$;
\item{(III)} the set of products $y_iz_j$, with $0\leq i\leq j-2\leq t-2$.
\end{list}
\noindent
We will denote by $\bar J_{s,t}$ the ideal associated with the matrix $\bar A_{s,t}$.\newline 
As shown in \cite{BM1}, Section 1, $J_{s,t}$ it is the defining ideal of a Cohen-Macaulay variety of minimal degree and it admits the prime decomposition 
$$J_{s,t}=J_0\cap J_1\cap\cdots\cap J_t,$$
\noindent where 
$$J_0=({\cal M}, {\cal D}_0),\quad\mbox{and }J_i=({\cal P}_i, {\cal D}_i)\quad\mbox{for }i=1,\dots, t,$$
\noindent with
$${\cal P}_i=\{x_1,\dots, x_s, y_0, \dots, y_{i-2}\}\quad\mbox{for }i=1,\dots, t,$$
$${\cal D}_i=\{z_{i+1}, \dots, z_t\}\quad\mbox{for }i=0,\dots, t.$$
\noindent
Thus the sequence of ideals $J_0, J_1, \dots, J_t$ fulfils condition 2 of Theorem 1 in \cite{M}, which implies that it is  {\it linearly joined}; this notion was introduced by Eisenbud, Green, Hulek and Popescu \cite{EGHP}, and was later intensively investigated by  Morales \cite{M}.  We also have 
\begin{equation}\label{height} {\rm height}\,J_{s,t}=s+t-1.
\end{equation}
\noindent
In the sequel, we will set $V_{s,t}=V(J_{s,t})$, and also $\bar V_{s,t}=V(\bar J_{s,t})$. Note  that $J_{s,1}$ has the same generators as $\bar J_{s,1}$, because the indeterminate $y_0$ does not appear in these generators. Consequently, we can identify $V_{s,1}$ with $\bar V_{s,1}$.  One should observe that, apart from this special case, for any  integers $s$ and $t$, $J_{s,t}$ does not denote a single ideal, but a class of ideals, namely the ideals attached to a matrix $A_{s,t}$ for some choice of the (identification between) the indeterminates  
$x_1,x_2,\dots, x_{2s},y_0, y_1, \dots, y_{t-1}$, $z_1, z_2,\dots, z_t$. The same remark applies to the variety $V_{s,t}$.\newline
For the proofs of the theorems on arithmetical ranks contained in the next section we will need the following two technical results, which are valid in any commutative unit ring $R$. 
\begin{lemma}\label{lemma1}{\rm (\cite{Ba3}, Corollary 3.2)} Let $\alpha_1,\alpha_2,\beta_1,\beta_2,\gamma\in R$. Then 
\begin{eqnarray*}&&\sqrt{(\alpha_1\beta_1-\alpha_2\beta_2,\ \beta_1\gamma,\ \beta_2\gamma)}=\\
&&\qquad\qquad\qquad\qquad\qquad\sqrt{(\alpha_1(\alpha_1\beta_1-\alpha_2\beta_2)+\beta_2\gamma,\ \alpha_2(\alpha_1\beta_1-\alpha_2\beta_2)+\beta_1\gamma)}.\end{eqnarray*}
\end{lemma}
The next claim is a slightly generalized version of \cite{Ba3}, Lemma 2.1 (which, in turn, extends \cite{SV}, Lemma, p. 249). The proof is the same as the one given in \cite{Ba3}, and will therefore be omitted here.
\begin{lemma}\label{lemma2} Let $P$ be a finite subset of elements of $R$, and  $I$ an ideal of $R$.  Let $P_1,\dots, P_r$ be subsets of $P$ such that
\begin{list}{}{}
\item[(i)] $\bigcup_{\ell=1}^rP_{\ell}=P$;
\item[(ii)] if $p$ and $p'$ are different elements of $P_{\ell}$ $(1\leq\ell\leq r)$ then $(pp')^m\in I+\left(\bigcup_{i=1}^{\ell-1}P_i\right)$ for some positive integer $m$. 
\end{list}
\noindent
Let $1\leq \ell\leq r$, and, for any $p\in P_{\ell}$, let $e(p)\geq1$ be an integer. We set $q_{\ell}=\sum_{p\in P_{\ell}}p^{e(p)}$. Then we get
$$\sqrt{I+(P)}=\sqrt{I+(q_1,\dots,q_r)}.$$
\noindent
where $(P)$ denotes the ideal of $R$ generated by $P$.
\end{lemma}
\section{The arithmetical rank for $s=2$: set-theoretic complete intersections}
In this section we will show that, for all $t\geq1$, the variety $V_{2,t}$ is a set-theoretic complete intersection. Recall that its defining ideal is  the ideal $J_{2,t}$ of $R=K[x_1, x_2, x_3, x_4, y_0, y_2, \dots, y_{t-1}, z_1, z_2,\dots, z_t]$, which is associated with the matrix 
$$A_{2,t}=\left(
\begin{array}{ccccccccccc}
x_1&x_2&\vline\ \vline&y_0&\vline\ \vline&y_1&\vline\ \vline&\cdots&\vline\ \vline&y_{t-1}\\
x_3&x_4&\vline\ \vline&z_1&\vline\ \vline&z_2&\vline\ \vline&\cdots&\vline\ \vline&z_t\\
\end{array}\right),
$$
\noindent
and is generated by the elements 
$$\begin{array}{ccccc}
x_1x_4-x_2x_3,& x_1z_1, &x_1z_2,&\dots,&x_1z_t,\\
&x_2z_1, &x_2z_2,&\dots,&x_2z_t,\\
&y_0z_2, &\dots,&\dots, &y_0z_t,\\
&&y_1z_3,&\dots, &y_1z_t,\\
&\dots,&\dots,&\dots, &y_{t-2}z_t.
\end{array}$$
\noindent
The next result generalizes Example 5 in \cite{Ba2}.
\begin{theorem}\label{theorem3} For all integers $t\geq 1$, {\rm ara}\,$J_{2,t}=t+1$, i.e., $V_{2,t}$ is a set-theoretic complete intersection.
\end{theorem}
\demo We proceed by induction on $t$, by showing that there are $F_1,\dots, F_{t+1}\in R= K[x_1, x_2, x_3, x_4, y_0, \dots, y_{t-1},z_1, z_2,\dots, z_t]$ such that
\begin{list}{}{}
\item{(a)} $\sqrt{(F_1,\dots, F_{t+1})}=J_{2,t}$,
\item{(b)} $F_1,F_2\in(x_1,x_2)$,
\item{(c)} $F_i\in(x_1, x_2, y_0, \dots, y_{i-3})$ for all $i=3,\dots, t+1$.
\end{list}
\noindent
For the induction basis consider the case where $t=1$. We have $J_{2,1}=(x_1x_4-x_2x_3, x_1z_1, x_2z_1)$. Set
\begin{equation}\label{F1F2}F_1=x_4(x_1x_4-x_2x_3)+x_2z_1,\qquad\qquad F_2=x_3(x_1x_4-x_2x_3)+x_1z_1.\end{equation}
\noindent
Then $F_1$ and $F_2$ fulfil condition (b) and, by virtue of Lemma \ref{lemma1}, they also fulfil condition (a). Now assume that $t\geq 2$ and suppose that $G_1,\dots, G_t$ are polynomials fulfilling the claim for $t-1$. By condition (b) we have  $G_1=Px_1-Qx_2$ for some $P,Q\in R$. Set
\begin{eqnarray*}
F_1&=&QG_1+x_1z_t\\
F_2&=&PG_1+x_2z_t\\
F_3&=&G_2+y_0z_t\\
\vdots&&\\
F_i&=&G_{i-1}+y_{i-3}z_t\\
\vdots&&\\
F_{t+1}&=&G_t+y_{t-2}z_t.
\end{eqnarray*}
\noindent
Then $F_1, F_2\in(G_1, x_1, x_2)\subset(x_1, x_2)$. Moreover, for all $i=3,\dots, t+1$,
$$F_i\in(G_{i-1}, y_{i-3})\subset (x_1,x_2, y_0,\dots, y_{i-4}, y_{i-3}),$$
\noindent
because $G_{i-1}$ fulfils condition (c). Hence $F_1,\dots, F_{t+1}$ fulfil conditions (b) and (c). Furthermore, by Lemma \ref{lemma1}, 
\begin{equation}\label{1}\sqrt{(F_1, F_2)}=\sqrt{(G_1, x_1z_t, x_2z_t)},\end{equation}
\noindent
and, for all $i=2,\dots, t$, the product of the two summands of $F_{i+1}$ is 
\begin{eqnarray*}
G_i\cdot y_{i-2}z_t\in(x_1, x_2, y_0, \dots, y_{i-3})\cdot(z_t)&=&(x_1z_t, x_2z_t)+(y_0z_t, \dots, y_{i-3}z_t),\\
&\subset&\sqrt{(F_1, F_2)}+(y_0z_t, \dots, y_{i-3}z_t),
\end{eqnarray*}
\noindent
where the first membership relation is true because $G_i$ fulfils condition (c). It follows that $(G_i\cdot y_{i-2}z_t)^m$ belongs to $(F_1, F_2)+(y_0z_t, \dots, y_{i-3}z_t)$ for some positive integer $m$.  Hence the assumption of Lemma \ref{lemma2} is fulfilled for $I=(F_1, F_2)$ and $P_i=\{G_{i+1}, y_{i-1}z_t\}$ ($i=1,\dots,t-1$).  Consequently, 
\begin{eqnarray*}
\sqrt{(F_1, F_2, F_3, \dots, F_{t+1})}&=&\sqrt{(F_1, F_2, G_2, \dots, G_t, y_0z_t, \dots, y_{t-2}z_t)}\\
&=&\sqrt{(G_1, G_2,\dots G_t, x_1z_t, x_2z_t, y_0z_t, \dots, y_{t-2}z_t)}\\
&=&J_{2, t-1}+(x_1z_t, x_2z_t, y_0z_t, \dots, y_{t-2}z_t)=J_{2,t},
\end{eqnarray*}
\noindent
where the second and the third equality follow from (\ref{1}) and induction, respectively. Thus $F_1,\dots, F_{t+1}$ fulfil condition (a) as well. This completes the proof.
\begin{remark}{\rm The polynomials $F_1,\dots, F_{t+1}$ defined in the proof of Theorem \ref{theorem3} still fulfil the required properties if in all  monomial summands $x_1z_t$, $x_2z_t$, $y_0z_t,\dots, y_{t-2}z_t$ the factors $z_t$ are raised to the same arbitrary positive power. This allows us, e.g., to replace the polynomials in (\ref{F1F2}) by 
$$F_1=x_4(x_1x_4-x_2x_3)+x_2z_1^2,\qquad   F_2=x_3(x_1x_4-x_2x_3)+x_1z_1^2,$$
\noindent 
which are homogeneous. Then, by a suitable adjustment of exponents, one can recursively construct a sequence of homogeneous polynomials $F_1,\dots, F_{t+1}$ for any $t\geq2$.}
\end{remark}
\begin{example}\label{Example1}{\rm Equalities (\ref{F1F2}) explicitly provide the defining polynomials for $V_{2,1}$. They are the starting point of the recursive procedure, described in the proof of Theorem \ref{theorem3}, which allows us to construct $t+1$ polynomials defining $V_{2,t}$, for any $t\geq2$. We perform the construction for $t=2,3$. First take $t=2$. We have
$$A_{2,2}=\left(\begin{array}{cccccc}
x_1&x_2&\vline\ \vline&y_0&\vline\ \vline& y_1\\
x_3&x_4&\vline\ \vline&z_1&\vline\ \vline& z_2
\end{array}\right),$$
and
$$J_{2,2}=(x_1x_4-x_2x_3, x_1z_1, x_1z_2, x_2z_1, x_2z_2, y_0z_2).$$
\noindent
Let us rewrite the polynomials given in (\ref{F1F2}): 
$$G_1=x_1x_4^2-x_2x_3x_4+x_2z_1,\qquad\qquad G_2=x_1x_3x_4-x_2x_3^2+x_1z_1.$$
\noindent
Then, with the notation of the proof of Theorem \ref{theorem3}, $P=x_4^2$ and $Q=x_3x_4-z_1$. Thus 
\begin{eqnarray*}
F_1&=&(x_3x_4-z_1)G_1+x_1z_2\\
&=&x_1x_3x_4^3-x_1x_4^2z_1-x_2x_3^2x_4^2+2x_2x_3x_4z_1-x_2z_1^2+x_1z_2,\\
F_2&=&x_4^2G_1+x_2z_2=x_1x_4^4-x_2x_3x_4^3+x_2x_4^2z_1+x_2z_2,\\
F_3&=&G_2+y_0z_2=x_1x_4x_3-x_2x_3^2+x_1z_1+y_0z_2
\end{eqnarray*}
are three defining polynomials for $V_{2,2}$. Now let $t=3$. We have
$$A_{2,3}=\left(\begin{array}{cccccccc}
x_1&x_2&\vline\ \vline&y_0&\vline\ \vline& y_1&\vline\ \vline& y_2\\
x_3&x_4&\vline\ \vline&z_1&\vline\ \vline& z_2&\vline\ \vline& z_3
\end{array}\right),$$
\noindent
and 
$$J_{2,3}=(x_1x_4-x_2x_3, x_1z_1, x_1z_2, x_1z_3, x_2z_1, x_2z_2, x_2z_3, y_0z_2, y_0z_3, y_1z_3).$$
\noindent
In order to obtain four defining polynomials for $V_{2,3}$ we take the above polynomials $F_1, F_2, F_3$ as $G_1, G_2, G_3$. Thus $P=x_3x_4^3-x_4^2z_1+z_2$ and $Q=x_3^2x_4^2-2x_3x_4z_1+z_1^2$. Hence, the four sought polynomials are
\begin{eqnarray*}
F_1&=&(x_3^2x_4^2-2x_3x_4z_1+z_1^2)G_1+x_1z_3=x_1x_3^3x_4^5-3x_1x_3^2x_4^4z_1+3x_1x_3x_4^3z_1^2\\
&&-x_1x_4^2z_1^3-x_2x_3^4x_4^4+4x_2x_3^3x_4^3z_1\\
&&-6x_2x_3^2x_4^2z_1^2+4x_2x_3x_4z_1^3-x_2z_1^4\\
&&+x_1x_3^2x_4^2z_2-2x_1x_3x_4z_1z_2+x_1z_1^2z_2+x_1z_3,\\
F_2&=&(x_3x_4^3-x_4^2z_1+z_2)G_1+x_2z_3=x_1x_3^2x_4^6-2x_1x_3x_4^5z_1+2x_1x_3x_4^3z_2\\
&&+x_1x_4^4z_1^2-2x_1x_4^2z_1z_2-x_2x_3^3x_4^5-x_2x_3^2x_4^2z_2\\
&&+3x_2x_3^2x_4^4z_1-3x_2x_3x_4^3z_1^2+2x_2x_3x_4z_1z_2\\
&&+x_2x_4^2z_1^3-x_2z_1^2z_2+x_1z_2^2+x_2z_3,\\
F_3&=&G_2+y_0z_3=x_1x_4^4-x_2x_3x_4^3+x_2x_4^2z_1+x_2z_2+y_0z_3,\\
F_4&=&G_3+y_1z_3=x_1x_3x_4-x_2x_3^2+x_1z_1+y_0z_2+y_1z_3.
\end{eqnarray*}}
\end{example}
\section{The arithmetical rank for $s\geq3$: upper and lower bounds}
The aim of this section is to show that, for $s\geq3$, the ideal $J_{s,t}$ is never a set-theoretic complete intersection. We will determine a lower bound for ara\,$J_{s,t}$, which shows that the difference between the arithmetical rank and the height strictly increases with $s$.  For our purpose we will need the following cohomological criterion by Newstead \cite{N}. 
\begin{lemma}\label{Newstead}{\rm (\cite{BS}, Lemma 3$^\prime$)} Let
$W\subset\tilde W$ be affine varieties. Let $d=\dim\tilde
W\setminus W$. If there are $r$ polynomials $F_1,\dots, F_r$ such
that $W=\tilde W\cap V(F_1,\dots,F_r)$, then 
$$\het^{d+i}(\tilde W\setminus W,{\bz}/m{\bz})=0\quad\mbox{ for all
}i\geq r$$ and for all $m\in{\bz}$ which are prime to {\rm char}\,$K$.
\end{lemma}
We refer to \cite{Mi0} or \cite{Mi} for the basic notions on \'etale cohomology. We are now ready to prove the first of the two main results of this section. 
\begin{theorem}\label{theorem5} For all integers $s\geq2$ and $t\geq1$
$${\rm ara}\,J_{s,t}\geq 2s+t-3.$$
\end{theorem}
\demo   For $s=2$ the claim is a trivial consequence of Theorem \ref{theorem3}. So let $s\geq 3$. It suffices to prove the claim for $\bar J_{s,t}$, because ara\,$J_{s,t}\geq$\,ara\,$\bar J_{s,t}$: in fact, given $r$ defining polynomials for $V_{s,t}$, they can be transformed in $r$ defining polynomials for $\bar V_{s,t}$ by performing on them the suitable identifications between the indeterminates. Let $p$ be a prime different from char\,$K$. According to Lemma \ref{Newstead}, it suffices to show that
\begin{equation}\label{2} \het^{4s+2t-4}(K^{2s+t}\setminus \bar V_{s,t}, \bz/p\bz)\ne0,\end{equation}
\noindent
since this will imply that $\bar V_{s,t}$ cannot be defined by $2s+t-4$ equations.
By Poincar\'e Duality (see \cite{Mi}, Theorem 14.7, p.~83) we have
\begin{equation}\label{3} {\rm Hom}_{\scriptstyle{\bz}/p{\bz}} (\het^{4s+2t-4}(K^{2s+t}\setminus \bar V_{s,t}, {\bz}/p{\bz}), {\bz}/p{\bz})\simeq  \hc^{4}(K^{2s+t}\setminus \bar V_{s,t}, {\bz}/p{\bz}),
\end{equation}
\noindent
where $\hc$ denotes \'etale cohomology with compact support. 
For the sake of simplicity, we will omit the coefficient group ${\bz}/p{\bz}$ henceforth. In view of (\ref{3}), it suffices to show that
\begin{equation}\label{3b} \hc^{4}(K^{2s+t}\setminus \bar V_{s,t})\ne0.\end{equation}
\noindent
 Let $W$ be the subvariety of $K^{2s+t}$ defined by the vanishing of $y_t$ and of all generators of $\bar J_{s,t}$ listed in Section 1 under (I) and (II), and those listed in (III) for which $j\leq t-1$. Then $W\subset \bar V_{s,t}$, and   
\begin{eqnarray}\label{6'}
&&\!\!\!\!\!\!\!\!\!\!\!\!\!\!\!\bar V_{s,t}\setminus W=\nonumber\\
&&\!\!\!\!\!\!\!\!\!\!\!\!\!\!\!\{(x_1, \dots x_{2s}, y_1,\dots, y_t)\vert x_1=\cdots=x_s=x_{2s}=y_1=\cdots=y_{t-2}=0, \ y_t\ne0\}\nonumber\\
&&\qquad\simeq K^s\times (K\setminus\{0\}),\end{eqnarray}
\noindent
It is well known that
\begin{equation}\label{4}\hc^i(K^r)\simeq\left\{\begin{array}{cl} {\bz}/p{\bz}&\mbox{if }i=2r\\
0&\mbox{else, }
\end{array}\right.
\end{equation}
\noindent
and 
\begin{equation}\label{5}\hc^i(K^r\setminus\{0\})\simeq\left\{\begin{array}{cl} {\bz}/p{\bz}&\mbox{if }i=1,2r\\
0&\mbox{else. }
\end{array}\right.
\end{equation}
\noindent
Moreover, in view of (\ref{6'}), by the K\"unneth formula  for \'etale cohomology (\cite{Mi}, Theorem 22.1), 
$$\hc^i(\bar V_{s,t}\setminus W)\simeq\displaystyle\bigoplus_{h+k=i}\hc^h(K^s)\otimes\hc^k(K\setminus\{0\}),$$
\noindent
so that, by (\ref{4}) and (\ref{5}), we have $\hc^i(\bar V_{s,t}\setminus W)\ne0$ if and only if $i=2s+1, 2s+2.$  But $4<2s\leq 2s+1$, so that, in particular
\begin{equation}\label{6}\hc^3(\bar V_{s,t}\setminus W)=\hc^4(\bar V_{s,t}\setminus W)=0.\end{equation}
We have a long exact sequence of \'etale cohomology with compact support (see \cite{Mi0}, Remark 1.30, p. 94):
$$\cdots\rightarrow\hc^3(\bar V_{s,t}\setminus W)\rightarrow \hc^4(K^{2s+t}\setminus \bar V_{s,t})\rightarrow \hc^4(K^{2s+t}\setminus W)\rightarrow 
\hc^4(\bar V_{s,t}\setminus W)\rightarrow\cdots.$$
\noindent
By (\ref{6}) it follows that
\begin{equation}\label{7} \hc^4(K^{2s+t}\setminus \bar V_{s,t})\simeq\hc^4(K^{2s+t}\setminus W).
\end{equation}
\noindent
Note that $W$ can be described as the variety of $K^{2s+t}$ defined by the vanishing of $y_t$ and of all polynomials defining $\bar V_{s,t-1}$ in $K^{2s+t-1}$.  Note that  a point of $K^{2s+t}$ belongs to $K^{2s+t}\setminus W$ if and only if it fulfils one of the two following complementary cases: 
\begin{list}{}{}
\item{-} either its $y_t$-coordinate is zero, and it does not annihilate all polynomials of $\bar J_{s,t-1}$, or
\item{-}  its $y_t$-coordinate is non zero.
\end{list} 
Therefore we have  
\begin{equation}\label{ast}K^{2s+t}\setminus W=(K^{2s+t-1}\setminus \bar V_{s, t-1})\cup Z,\end{equation}
 where the union is disjoint, and $Z$ is the open subset given by
\begin{equation}\label{decomposition}Z=K^{2s+t-1}\times (K\setminus\{0\}).\end{equation}
\noindent We thus have a long exact sequence of \'etale cohomology with compact support:
\begin{equation}\label{9}\cdots\rightarrow\hc^4(Z)\rightarrow \hc^4(K^{2s+t}\setminus W)\rightarrow \hc^4(K^{2s+t-1}\setminus \bar V_{s, t-1})\rightarrow \hc^5(Z)\rightarrow\cdots.\end{equation}
\noindent
By  the K\"unneth formula for \'etale cohomology, (\ref{4}), (\ref{5}) and (\ref{decomposition}), we have $\hc^i(Z)\ne 0$ if and only if $i=4s+2t-1, 4s+2t$. But $4s+2t-1>5$, whence, in particular,
$$\hc^4(Z)=\hc^5(Z)=0.$$
\noindent
 It follows that (\ref{9}) gives rise to an isomorphism: 
$$\hc^4(K^{2s+t}\setminus W)\simeq \hc^4(K^{2s+t-1}\setminus \bar V_{s, t-1}).$$
\noindent
Hence, in view  of (\ref{7}), claim (\ref{3b}) follows by induction on $t$ if it is true that 
\begin{equation}\label{basis}\hc^4(K^{2s}\setminus V_{s, 0})\ne 0,\end{equation}
\noindent
where $V_{s,0}\subset K^{2s}$ denotes the variety defined by the vanishing of the two-minors of the first big block of $A_{s,t}$. But according to \cite{BS}, Lemma 2$'$, $\het^{4s-4}(K^{2s}\setminus V_{s,0})\ne 0$, from which (\ref{basis}) can be deduced by Poincar\'e Duality. This completes the proof of the theorem.
\begin{remark}{\rm According to (\ref{height}) and Theorem \ref{theorem5}, the difference between the arithmetical rank and the height of $J_{s,t}$ is at least $2s+t-3-(s+t-1)=s-2$. Thus it strictly increases with $s$. In view of Theorem \ref{theorem3}, it is zero if and only if $s=2$. 
}
\end{remark}
\begin{corollary} The variety $V_{s,t}$ is a set-theoretic complete intersection if and only if $s=2$.
\end{corollary}
Next we give an upper bound for ara\,$J_{s,t}$. In the sequel, for the sake of simplicity, we will denote by $[ij]$ $(1\leq i<j\leq s)$ the minor formed by the $i$th and the $j$th column of $A_{s,t}$. We will call $I_s$ the ideal generated by these minors (it is the defining ideal of the variety $V_{s,0}$ mentioned in the proof of Theorem \ref{theorem5}). Moreover, for all $k=1,\dots, 2s-3$, we set 
$$S_k=\sum_{i+j=k+2}[ij].$$
\noindent
 We preliminarily recall an important result by Bruns et al.
\begin{theorem}\label{theorem6}{\rm (\cite{BS}, Theorem 2 and \cite{BV}, Corollary 5.21)} With the notation just introduced,
$${\rm ara}\,I_s=2s-3,$$
\noindent
and
$$I_s=\sqrt{(S_1,\dots, S_{2s-3})}.$$ 
\end{theorem}
We can now prove the second result of this section.
\begin{theorem}\label{theorem7} For all integers $s\geq2$ and $t\geq1$,
$${\rm ara}\,J_{s,t}\leq 2s+t-2.$$
\end{theorem}
\demo Again, in view of Theorem \ref{theorem3}, it suffices to prove the claim for $s\geq 3$. Let $L_{s,t}$ be the ideal generated by the products listed in Section 1 under (II) and (III). For convenience of notation we set
\begin{eqnarray*} \xi_i&=&x_i\qquad\qquad(1\leq i\leq s)\\
\xi_i&=&y_{i-s-1}\qquad(s+1\leq i\leq s+t).
\end{eqnarray*}
\noindent
In other words, the entries of the first row of $A_{s,t}$ are denoted by $\xi_1,\dots, \xi_{s+t}$, and the monomial generators of $L_{s,t}$ are
\begin{equation}\label{Lst}\xi_iz_j,\qquad\mbox{where }1\leq i\leq s+t-1,\  i-s+1\leq j\leq t.\end{equation}
\noindent
Let 
$$T_h=\sum_{i=1}^{s+t-1}\xi_iz_{i+t-h}\qquad\qquad (1\leq h\leq s+t-1),$$
\noindent
where we have set  $z_j=0$ for $j\not\in\{1,\dots, t\}$. Then the set of non zero monomial summands  in $T_1,\dots, T_{s+t-1}$ coincides with the set of monomial generators of $L_{s,t}$, as the following elementary argument shows. On the one hand, given a non zero monomial summand $\xi_iz_{i+t-h}$ of some $T_h$, it holds
$$i-s+1=i+t-s-t+1\leq i+t-h,$$
\noindent
so that  $\xi_iz_{i+t-h}$ is of the form (\ref{Lst}). On the other hand, given a monomial $\xi_iz_j$ as in (\ref{Lst}), we have $j=i+t-h$ for $h=i+t-j$, where $j\leq t$ and $i-s+1\leq j$. Therefore,
$$1\leq i\leq h\leq i+t-(i-s+1)=s+t-1,$$
\noindent
which implies that $\xi_iz_j$ is a monomial summand of $T_h$.\newline
Moreover, $T_1=\xi_1z_t$. Now consider, for any $h$ such that $1\leq h\leq s+t-1$,   the product  of two non zero distinct monomial summands of $T_h$: it is of the form $\xi_pz_{p+t-h}\xi_qz_{q+t-h}$ with $1\leq p<q\leq s+t-1$. Hence it is divisible by $\xi_pz_{q+t-h}=\xi_pz_{p+t-(h+p-q)}$, which is one of the non zero monomial summands of $T_{h+p-q}$.   Since $q+t-h\leq t$, we have $h-q\geq0$, whence it follows that $1\leq p\leq h+p-q< h$. Thus the assumption of Lemma \ref{lemma2} is fulfilled if we take $I=(T_1)$, $P_h$ equal to the set of all non zero monomial summands of $T_h$ and $q_h=T_h$ for $h=2,\dots,s+t-1$. Therefore
\begin{equation}\label{eq3}L_{s,t}=\sqrt{(T_1,\dots, T_{s+t-1})}.\end{equation}
\noindent
For some arbitrarily fixed $\ell$ with $1\leq \ell\leq 2s-3$, let $[ij]$ be a summand of $S_{\ell}$. Then the monomial terms of $[ij]$ are of the form
\begin{equation}\label{eq4} \xi_ux_v,\qquad\mbox{where $1\leq u\leq \ell+1$.}\end{equation}
\noindent
For some fixed $h$ such that $1\leq h\leq s+t-1$, let $\xi_iz_{i+t-h}$ be a non zero monomial summand of $T_h$. Then $i+t-h\geq 1$ implies that 
\begin{equation}\label{eq5} h-i\leq t-1. \end{equation}
\noindent
For all $\ell=1,\dots, s-2$ let
\begin{equation}\label{Ul}U_{\ell}=S_{\ell}+T_{\ell+t+1}.\end{equation}
\noindent
Then, if $\xi_ux_v$ is a monomial term in $S_{\ell}$ and $\xi_iz_{i+t-(\ell+t+1)}$ a non zero monomial summand in $T_{\ell +t+1}$, their product is divisible by
\begin{equation}\label{T}\xi_uz_{i+t-(\ell +t+1)}=\xi_uz_{u+t-(\ell+t+1+u-i)}.\end{equation}
\noindent
Set $h'=\ell+t+1+u-i.$ Now, according to (\ref{eq4}), $u\leq\ell +1$,  so that, applying (\ref{eq5}) for $h=\ell+t+1$, we obtain $h'=\ell+t+1-i+u\leq t-1+\ell+1=\ell+t$. On the other hand, since $z_{i+t-(\ell+t+1)}\ne 0$, we have  $i+t-(\ell+t+1)\leq t$, i.e., $\ell+t+1\geq i$.  This implies that $h'=\ell+t+1+u-i\geq u\geq 1.$ Thus (\ref{T}) shows that the product of each two-minor appearing as a summand in $S_{\ell}$ and each non zero monomial summand of $T_{\ell+t+1}$ is divisible by a monomial summand of $T_{h'}$, for some $h'$ such that $1\leq h'< \ell +t+1$. Thus Lemma \ref{lemma1} 
applies to $I=(T_1,\dots, T_{t+1})$, $P'_{\ell}=\{S_{\ell}, T_{\ell+t+1}\}$  and $q'_{\ell}=U_{\ell}$ for $\ell=1,\dots, s-2$, whence, in view of (\ref{Ul}), we conclude that
\begin{eqnarray*}
\sqrt{(T_1,\dots, T_{t+1}, U_1,\dots, U_{s-2})}&=&\\
\sqrt{(T_1,\dots, T_{t+1}, T_{t+2}, \dots, T_{s+t-1}, S_1,\dots, S_{s-2})}&=&\sqrt{L_{s,t}+(S_1, \dots, S_{s-2})},
\end{eqnarray*}
\noindent
where the last equality is a consequence of (\ref{eq3}).
Thus we have 
\begin{eqnarray}\label{ideal}
\sqrt{(T_1,\dots, T_{t+1}, U_1,\dots, U_{s-2}, S_{s-1}, \dots, S_{2s-3})}&=&\\
\sqrt{L_{s,t}+(S_1, \dots, S_{2s-3})}&=&\sqrt{L_{s,t}+I_s}=J_{s,t},\nonumber
\end{eqnarray}
\noindent
where the second equality follows from Theorem \ref{theorem6}.
Since  the number of generators of the ideal  in (\ref{ideal}) is $t+1+2s-3=2s+t-2$, this completes the proof.
\par\medskip\noindent
The gap between the lower bound given in Theorem \ref{theorem5} and the upper bound given in Theorem \ref{theorem7} is equal to 1. Theorem \ref{theorem3} also shows that the lower bound is sharp.
\begin{corollary}\label{tight} For all integers $s\geq 2$ and $t\geq 1$,
$$2s+t-3\leq\,{\rm ara}\,J_{s,t}\leq 2s+t-2.$$
\noindent
If $s=2$, then the first inequality is an equality.
\end{corollary}
There are other cases where the lower bound is sharp. In fact it is the exact value of ara\,$J_{s,1}$ for $s=3,4,5$, i.e., we have ara\,$J_{3,1}=4$, ara\,$J_{4,1}=6$, ara\,$J_{5,1}=8$. This is what we are going to show in the next example: it will suffice to produce, in the three aforementioned cases, 4, 5 and 6 defining polynomials, respectively. 
\begin{example}{\rm With the notation introduced above, we have 
$$A_{3,1}=\left(\begin{array}{ccccc}
x_1&x_2&x_3&\vline\ \vline& y_0\\
x_4&x_5&x_6&\vline\ \vline& z_1
\end{array}
\right),$$
\noindent and
$$J_{3,1}=([12], [13], [23], x_1z_1, x_2z_1, x_3z_1),$$
where
$$[12]=x_1x_5-x_2x_4,\qquad [23]=x_2x_6-x_3x_5,\qquad [13]=x_1x_6-x_3x_4.$$
\noindent 
We show that four defining polynomials are:
\begin{eqnarray*}
F_1&=&[23]\\
F_2&=&\hphantom{[13]+}\,\,x_1z_1+x_4[12]\\
F_3&=&[13]+x_2z_1+x_5[12]\\
F_4&=&\hphantom{[13]+}\,\,x_3z_1+x_6[12]
\end{eqnarray*}
Since $F_1, F_2, F_3, F_4\in J_{3,1}$, by virtue of Hilbert's Nullstellensatz it suffices to prove that every ${\bf v}=(x_1, \dots, x_6, z_1)\in K^7$ which annihilates all four polynomials annihilates all generators of $J_{3,1}$. In the sequel, we will use, when this does not cause any confusion, the same notation for the polynomials and for their values at ${\bf v}$. From $F_1=0$ we immediately get $[23]=0$. Moreover, since ${\bf v}$ annihilates $F_2, F_3, F_4$, we have that the triple $([13], z_1, [12])$ is a solution of the $3\times 3$ system of homogeneous linear equations associated with the matrix 
$$\left(\begin{array}{ccc}
0&x_1&x_4\\
1&x_2&x_5\\
0&x_3&x_6
\end{array}
\right),
$$
whose determinant is $$\Delta=-x_1x_6+x_3x_4=-[13].$$ By Cramer's Rule, whenever $\Delta\ne0$, the only solution is the trivial one, so that, in particular, $[13]=0$, a contradiction. Thus we always have  $\Delta=0$, i.e., $[13]=0$. Hence, in view of Lemma \ref{lemma1},  $F_2=F_3=0$ implies that $[12]=x_1z_1=x_2z_1=0$. Consequently, $F_4=0$ implies that $x_3z_1=0$. Thus ${\bf v}$ annihilates all generators of $J_{3,1}$, as required. This shows that ara\,$J_{3,1}=4$. \newline
Now consider
$$A_{4,1}=\left(\begin{array}{cccccc}
x_1&x_2&x_3&x_4&\vline\ \vline& y_0\\
x_5&x_6&x_7&x_8&\vline\ \vline& z_1
\end{array}
\right).$$
By Theorem \ref{theorem6} we have
\begin{eqnarray}\label{radical}
J_{4,1}&=&([12], [13], [14], [23], [24], [34], x_1z_1, x_2z_1, x_3z_1, x_4z_1),\nonumber\\
&=&\sqrt{([12], [13], [14]+[23], [24], [34], x_1z_1, x_2z_1, x_3z_1, x_4z_1)},
\end{eqnarray}
where
$$[12]=x_1x_6-x_2x_5,\qquad [13]=x_1x_7-x_3x_5,\qquad [14]=x_1x_8-x_4x_5,$$
$$[23]=x_2x_7-x_3x_6,\qquad [24]=x_2x_8-x_4x_6,\qquad [34]=x_3x_8-x_4x_7.$$
\noindent 
Six defining polynomials are:
\begin{eqnarray*}
F_1&=&[24]\\
F_2&=&[14]+[23]\\
F_3&=&[34]+\hphantom{[13]+}\,\,\,x_1z_1+x_5[12]\\
F_4&=&\hphantom{[34]+}\,\,[13]+\,x_2z_1+x_6[12]\\
F_5&=&\hphantom{[34]+[13]+}\,\,\,x_3z_1+x_7[12]\\
F_6&=&\hphantom{[34]+[13]+}\,\,\,x_4z_1+x_8[12].
\end{eqnarray*}
\noindent
Suppose that all these polynomials vanish at ${\bf v}=(x_1, \dots, x_8, z_1)\in K^9$. We show that then ${\bf v}$ annihilates all generators of the ideal appearing under the radical sign in (\ref{radical}). From $F_1=F_2=0$ we get that $[24]=[14]+[23]=0$.  Moreover, since ${\bf v}$ annihilates $F_3,\dots, F_6$, we have that the quadruple $([34], [13], z_1, [12])$ is a solution of the $4\times 4$ system of homogeneous linear equations associated with the matrix 
$$\left(\begin{array}{cccc}
1&0&x_1&x_5\\
0&1&x_2&x_6\\
0&0&x_3&x_7\\
0&0&x_4&x_8
\end{array}
\right),
$$
whose determinant is $$\Delta=x_3x_8-x_4x_7=[34].$$ By Cramer's Rule, if $\Delta\ne0$, the only solution is the trivial one, so that, in particular, $[34]=0$, a contradiction. Hence we always have $\Delta=0$, i.e., $[34]=0$.  Hence, in analogy to what has been shown for $J_{3,1}$,  $F_3=F_4=F_5=0$ implies that $[13]=x_1z_1=x_2z_1=x_3z_1=[12]=0$. Consequently, $F_6=0$ implies that $x_4z_1=0$. Thus ${\bf v}$ annihilates all generators of the ideal in (\ref{radical}), as required. This shows that ara\,$J_{4,1}=6$. \newline
Finally consider
$$A_{5,1}=\left(\begin{array}{ccccccc}
x_1&x_2&x_3&x_4&x_5&\vline\ \vline& y_0\\
x_6&x_7&x_8&x_9&x_{10}&\vline\ \vline& z_1
\end{array}
\right).$$
By Theorem \ref{theorem6} we have
\begin{eqnarray}\label{radical2}
&&\!\!\!\!\!J_{5,1}=\nonumber\\
&&\!\!\!\!\!\!\!([12], [13], [14], [15], [23], [24], [25], [34], [35], [45], x_1z_1, x_2z_1, x_3z_1, x_4z_1, x_5z_1)=\nonumber\\
&&\!\!\!\!\!\!\!\!\!\!\!\!\sqrt{([12], [13], [14]+[23], [15]+[24], [25]+[34], [35], [45], x_1z_1, x_2z_1, x_3z_1, x_4z_1, x_5z_1)}\nonumber.\\
&&
\end{eqnarray}
where
$$[12]=x_1x_7-x_2x_6,\ [13]=x_1x_8-x_3x_6,\ [14]=x_1x_9-x_4x_6,\ [15]=x_1x_{10}-x_5x_6,$$
$$[23]=x_2x_8-x_3x_7,\ [24]=x_2x_9-x_4x_7,\ [25]=x_2x_{10}-x_5x_7,\ [34]=x_3x_9-x_4x_8,$$
$$[35]=x_3x_{10}-x_5x_8,\  [45]=x_4x_{10}-x_5x_9.$$
\noindent 
Eight defining polynomials are:
\begin{eqnarray*}
F_1&=&[14]+[23]\\
F_2&=&[15]+[24]\\
F_3&=&[25]+[34]\\
F_4&=&\hphantom{[45]+}\,\,[35]\,\hphantom{[35]}\,\,\,\,+\,x_1z_1+x_6[12]\\
F_5&=&\hphantom{[45]+[13]+}[13]+\,x_2z_1+x_7[12]\\
F_6&=&[45]\hphantom{[13]+[35]+}\,+x_3z_1+x_8[12]\\
F_7&=&\hphantom{[45]+[13]+[35]+}\,\,x_4z_1+x_9[12]\\
F_8&=&\hphantom{[45]+[13]+[35]+}\,x_5z_1+x_{10}[12]\\
\end{eqnarray*}
\noindent
Suppose that all these polynomials vanish at ${\bf v}=(x_1, \dots, x_{10}, z_1)\in K^{11}$. We show that then ${\bf v}$ annihilates all generators of the ideal appearing under the radical sign in (\ref{radical2}). From $F_1=F_2=F_3=0$ we get that 
$[14]+[23]=[15]+[24]=[25]+[34]=0$.  Moreover, since ${\bf v}$ annihilates $F_4,\dots, F_8$, we have that the 5-uple $([45], [35], [13], z_1, [12])$ is a solution of the  $5\times 5$ system of homogeneous linear equations associated with the matrix  
$$\left(\begin{array}{ccccc}
0&1&0&x_1&x_6\\
0&0&1&x_2&x_7\\
1&0&0&x_3&x_8\\
0&0&0&x_4&x_9\\
0&0&0&x_5&x_{10}
\end{array}
\right),
$$
whose determinant is $$\Delta=x_4x_{10}-x_5x_9=[45].$$ By Cramer's Rule, if $\Delta\ne0$, the only solution is the trivial one, so that, in particular, $[45]=0$, a contradiction. Hence we always have $\Delta=0$, i.e., $[45]=0$.  Hence, in analogy to what has been shown for $J_{4,1}$,  $F_4=F_5=F_6=F_8=0$ implies that $[13]=[35]=x_1z_1=x_2z_1=x_3z_1=x_5z_1=[12]=0$. Consequently, $F_7=0$ implies that $x_4z_1=0$. Thus ${\bf v}$ annihilates all generators of the ideal in (\ref{radical2}), as required. This shows that ara\,$J_{5,1}=8$.}
\end{example}
\section{On cohomological dimensions}
 Recall that, for any proper ideal $I$ of $R$, the {\it (local) cohomological dimension} of $I$ is defined as the number
\begin{eqnarray*}{\rm cd}\,I&=&\max\{i\vert H^i_I(R)\ne0\},\\
&=&\min\{i\vert H^j_I(M)=0\mbox{ for all }j>i\mbox{ and all $R$-modules $M$}\}.
\end{eqnarray*}
\noindent 
where $H^i_I$ denotes the $i$th right derived functor of the local cohomology functor $\Gamma_I$;  we refer to  Brodmann and Sharp \cite{BrSh} or to Huneke and Taylor \cite{Hu} for an extensive exposition of this subject. In this section we will determine ${\rm cd}\,J_{s,t}$ for all integers $s\geq2$ and $t\geq1$. We will use the following technical results on De Rham   ($H_{\rm DR})$ and singular cohomology ($H$) with respect to the coefficient field $\bc$. The first involves sheaf cohomology (see \cite{BrSh}, Chapter 20, or \cite{Hu}, Section 2.3) with respect to the structure sheaf $\tilde R$ of $K^N$.  The second result is analogous to Lemma \ref{Newstead}.
\begin{lemma}\label{lemma8}{\rm (\cite{Hart2}, Proposition 7.2)} Let $V\subset K^N$ be a non singular complex variety of dimension $d$ such that $H^i(V,\tilde R)=0$ for all $i\geq r$. Then $H_{\rm DR}^i(V, \bc)=0$ for all $i\geq d+r$. 
\end{lemma}
\begin{lemma}\label{Newstead2}{\rm (\cite{BS}, Lemma 3)} Let
$W\subset\tilde W$ be affine complex  varieties such that $\tilde W\setminus W$ is non singular of pure dimension $d$. If there are $r$ polynomials $F_1,\dots, F_r$ such
that $W=\tilde W\cap V(F_1,\dots,F_r)$, then 
$$H^{d+i}(\tilde W\setminus W,\bc)=0$$
\noindent
for all $i\geq r$.
\end{lemma}
We also recall that, for every proper ideal $I$ of $R$,
\begin{equation}\label{cdara}{\rm cd}\,I\leq\,{\rm ara}\,I,\end{equation}
which is shown in \cite{Hart1}, Example 2, p. 414 (and also in \cite{BrSh}, Corollary 3.3.3, and in \cite{Hu}, Theorem 4.4). Equality holds if $I$ is generated by a regular sequence, in which case the aritmetical rank is equal to the length of that sequence.\newline
In the proof of the next result we will use the well-known characterization of local cohomology in terms of Koszul (or \u{C}ech) cohomology (see \cite{BrSh}, Section 5.2, or \cite{Hu}, Section 2.1). Let $u_1, \dots, u_h\in R$ be non zero generators of the proper ideal $I$ of $R$.   For all $S\subset\{1,\dots, h\}$ let $R_S$ denote the localization of $R$ with respect to the multiplicative set of $R$ generated by $\{u_i\vert i\in S\}$; set $R_{\emptyset}=R$. Then, according to \cite{Hu}, Theorem 2.10, or \cite{BrSh}, Theorem 5.1.19, for all $i\geq 0$, $H^i_I(R)$ is isomorphic to the $i$th cohomology module of a cochain complex $(C^{\cdot}, \phi_{\cdot})$ of $R$-modules constructed as follows (see \cite{BrSh}, Proposition 5.1.5).   For all $i\geq 0$, set 
$$C^i=\displaystyle\bigoplus_{S\subset\{1,\dots, h\}\atop{\vert S\vert=i}}R_S.$$
\noindent
Given any $\alpha\in C^i$, for all $i\geq1$ and all $S\subset\{1,\dots, h\}$ such that $\vert S\vert=i$,  $\alpha_S$ will denote the component of $\alpha$ in $R_S$. The map $\phi_{i-1}:C^{i-1}\to C^i$ is defined in such a way that, for every $\alpha\in C^{i-1}$, and for all $S\subset\{1,\dots, h\}$ for which $\vert S\vert=i$, 
$$\phi_{i-1}(\alpha)_S=\sum_{k\in S}c_{S,k}\frac{\alpha_{S\setminus\{k\}}}{1},$$
where $c_{S,k}\in\{-1, 1\}$ and $\frac{\alpha_{S\setminus\{k\}}}{1}$ is the image of $\alpha_{S\setminus\{k\}}$ under the localization map $R_{S\setminus\{k\}}\to R_S$. 

\begin{lemma}\label{regular} Let $z$ be one of the indeterminates of $R$ and let $I$ be an ideal of $R$ generated by polynomials in which $z$ does not occur. Then, for all $i\geq 0$,
\begin{list}{}{}
\item{(i)}  $z$ is regular on $H^i_I(R)$;
\item{(ii)} if $H^i_I(R)\ne 0$, then $H^i_I(R)\ne zH^i_I(R)$.
\end{list}
\end{lemma}
\demo Let $u_1,\dots, u_h$ be non zero generators of $I$ not containing the indeterminate $z$. Let $S\subset\{1,\dots, h\}$. In this proof,  we will say that an element $a\in R_S$ does not contain the indeterminate $z$ if 
$$a=\frac{f}{\displaystyle\prod_{k\in S}u_k^{s_k}},$$
\noindent
where $f\in R$ is a polynomial not containing the indeterminate $z$. This definition is of course independent of the choice of $f$ and of the exponents $s_k$.  Moreover, there is a unique decomposition
$$a=\bar a +z\tilde a$$
such that $\bar a,\tilde a\in R_S$ and $\bar a$ does not contain $z$. Given $\alpha\in C_i$, for some $i\geq 0$, we will set $\bar\alpha=(\bar\alpha_S)_S$ and $\tilde\alpha=(\tilde\alpha_S)_S$, so that we have 
\begin{equation}\label{zdecomposition}
\alpha=\bar\alpha+z\tilde\alpha.
\end{equation}
\noindent
We will say that $\alpha$ is $z$-free whenever $\alpha=\bar\alpha$. The decomposition (\ref{zdecomposition}) is unique, and will be called the $z$-decomposition of $\alpha$. From the definition of $\phi_i$ it immediately follows that if $\alpha$ is $z$-free, so is $\phi_i(\alpha)$. Hence 
\begin{equation}\label{zdecompositionphi}
\phi_i(\alpha)=\phi_i(\bar\alpha)+z\phi_i(\tilde\alpha)
\end{equation}
is the $z$-decomposition of $\phi_i(\alpha)$. We thus have, for all $\alpha\in C_i$, 
\begin{equation}\label{28}\alpha\in\,{\rm Ker}\,\phi_i\Longleftrightarrow \bar\alpha,\tilde\alpha\in\,{\rm Ker}\,\phi_i,\end{equation}
\begin{equation}\label{28'}\alpha\in\,{\rm Im}\,\phi_{i-1}\Longleftrightarrow \bar\alpha,\tilde\alpha\in\,{\rm Im}\,\phi_{i-1}.\end{equation}
\noindent
Let $\alpha\in C_i$. First suppose that $z\alpha\in\,{\rm Im}\,\phi_{i-1}$. Then, for some $\beta\in C_{i-1}$,    $z\alpha=\phi_{i-1}(\beta)=\phi_{i-1}(\bar \beta)+z\phi_{i-1}(\tilde\beta)$, whence $\phi_{i-1}(\bar\beta)=0$ and $\alpha=\phi_{i-1}(\tilde\beta)$. Thus $\alpha\in{\rm Im}\,\phi_{i-1}.$ This proves part (i) of the claim.  Now suppose that  $H_I^i(R)\ne0$. Then there is $\alpha\in\,{\rm Ker}\,\phi_i$ such that $\alpha\not\in\,{\rm Im}\,\phi_{i-1}$. From (\ref{28}) and (\ref{28'}) we can easily deduce that one can choose $\alpha$ to be $z$-free. Suppose that $\alpha\in\,{\rm Im}\,\phi_{i-1}+z\,{\rm Ker}\,\phi_i$, i.e., $\alpha=\phi_{i-1}(\beta)+z\alpha'$ for some $\beta\in C_{i-1}$, $\alpha\in\,{\rm Ker}\,\phi_i$. By the uniqueness of the $z$-decomposition of $\alpha$ it follows that $\alpha=\phi_{i-1}(\bar\beta)$, a contradiction. This shows that  ${\rm Ker}\,\phi_i\neq\,{\rm Im}\,\phi_{i-1}+z\,{\rm Ker}\,\phi_i$, so that $H_I^i(R)\ne zH_I^i(R).$ This shows part (ii) of the claim and completes the proof.
\par\medskip\noindent
\begin{lemma}\label{plusone} Let $I$ be a proper ideal of $R$ generated by polynomials in which the indeterminate $z$ does not occur. Then
$${\rm cd}\,(I+(z))={\rm cd}\,I+1.$$
\end{lemma}
\demo  The claim for $I=(0)$ is true because, by the observation following (\ref{cdara}), we have that ${\rm cd}\,(z)=1$. So assume that $I\ne(0)$. 
Set $d={\rm cd}\,I$. We prove the claim by showing the two inequalities separately. We have the following exact sequence, the so-called Brodmann sequence (see \cite{Hu}, Theorem 3.2):
$$\cdots\to H_I^{i-1}(R_z)\to H_{I+(x)}^i(R)\to H_I^i(R)\to H_{I}^i(R_z)\to\cdots\,$$
\noindent
We deduce that $H_{I+(x)}^i(R)=0$ whenever $H_I^{i-1}(R_z)=H_I^i(R)=0$, which is certainly the case if $i>d+1$. It follows that ${\rm cd}\,(I+(z))\leq\,d+1$.
 By virtue of Lemma \ref{regular}, part (i), multiplication by $z$ on $H_I^d(R)$ gives rise to a short exact sequence
$$0\to H_I^d(R)\to H_I^d(R)\to H_I^i(R)/zH_I^d(R)\to0$$
from which, in turn, we obtain the long exact sequence of local cohomology:
\begin{equation}\label{longcohomology}\cdots\to H_{(z)}^0(H_I^d(R))\to H_{(z)}^0(H_I^d(R)/zH_I^d(R))\to H_{(z)}^1(H_I^d(R))\to\cdots.\end{equation}
\noindent
Now $H_{(z)}^0(H_I^d(R))\simeq \Gamma_{(z)}(H_I^d(R))=0$, because $z$ is regular on $H_I^d(R)$ by Lemma \ref{regular}, part (i). Moreover, 
$H_{(z)}^0(H_I^d(R)/z H_I^d(R))\simeq \Gamma_{(z)}(H_I^d(R)/z H_I^d(R))=H_I^d(R)/z H_I^d(R)$, since $H_I^d(R)/z H_I^d(R)$ is annihilated by $z$. Hence, by Lemma \ref{regular}, part (ii), we deduce that $H_{(z)}^0(H_I^d(R)/z H_I^d(R))\neq0$. Therefore, from (\ref{longcohomology}) it follows that 
$$H_{(z)}^1(H_I^d(R))\ne 0,$$
\noindent whereas from (\ref{cdara}) we know that
$$H_{(z)}^i(H_I^d(R))=0 \qquad\qquad\mbox{for all }i>1.$$
We have a Grothendieck spectral sequence for local cohomology (see \cite{R}, Theorem 11.38, or \cite{MC}, Theorem 12.10), 
$$E_2^{pq}=H_{(z)}^p(H_I^q(R))\Rightarrow H_{I+(z)}^{p+q}(R).$$
\noindent
The maximum value of $p+q$ for which $E_2^{pq}\ne0$ is $d+1$ and is obtained only for $p=1$ and $q=d$. 
Thus we get
$$H_{I+(z)}^{d+1}(R)\ne0,$$
\noindent
which yields ${\rm cd}\,(I+(z))\geq d+1$. This completes the proof.
\par\medskip\noindent
Before coming to the main result of this section, we first show one special case of its claim. This case deserves to be considered separately, because it is the only one where the cohomological dimension is independent of the characteristic of the ground field. The next proposition is an application of a recent result by Morales \cite{M}.
\begin{proposition}\label{special} Let $t\geq1$ be an integer. Then
 $${\rm cd}\,J_{2,t}=t+1.$$
\end{proposition}  
\demo
We refer to the prime decomposition of $J_{2,t}$ given in Section 1. Since ${\cal M}=\{x_1x_4-x_2x_3\}$, all ideals $J_0, J_1, \dots, J_t$ are complete intersections. According to \cite{M}, Theorem 4, this implies that 
\begin{equation}\label{morales}{\rm cd}\,J_{2,t}=\max_{j=1,\dots, t}\dim_K(\langle {\cal P}_i\rangle+\langle {\cal D}_{i-1}\rangle)-1.\end{equation}
\noindent Here the angle brackets denote linear spaces. It is evident from their definition that, for all $i=1,\dots, t$, ${\cal P}_i$ and ${\cal D}_{i-1}$ are disjoint sets of  $i+1$ and $t-i+1$ indeterminates, respectively. Hence $\dim_K(\langle {\cal P}_i\rangle+\langle {\cal D}_{i-1}\rangle)=\vert {\cal P}_i\vert +\vert {\cal D}_{i-1}\vert=t+2$ for all $i=1,\dots, t$, whence, in view of (\ref{morales}), the claim follows. 
\par\medskip\noindent
\begin{theorem}\label{theorem9} Let $s\geq2$ and $t\geq1$ be integers. Then
\begin{list}{}{}
\item{(a)} if char\,$K>0$, {\rm cd}\,$J_{s,t}=s+t-1$,
\item{(b)} if char\,$K=0$, {\rm cd}\,$J_{s,t}=2s+t-3$.
\end{list}
\end{theorem}
\demo Claim (a) follows from (\ref{height}) and \cite{PS}, Proposition 4.1, p. 110, since $J_{s,t}$ is Cohen-Macaulay. We prove claim (b) by  induction on $t$. Suppose that char\,$K=0$. The claim for $s=2$  and any integer $t\geq1$ is given by Proposition \ref{special}.  Next we consider the case where $s=3$ and $t=1$. We have  cd\,$J_{3,1}\leq 4$: this follows from (\ref{cdara}), since we have seen in Example \ref{Example1} that ara\,$J_{3,1}=4$. The same inequality has also been proven, by other means, in \cite{Ba2}, Example 6. In order to prove the opposite inequality, we have to show that
\begin{equation}\label{claimc}H_{J_{3,1}}^4(R)\ne0.\end{equation}
By virtue of the flat basis change property of local cohomology (see \cite{BrSh}, Theorem 4.3.2, or \cite{Hu}, Proposition 2.11 (1)), if this is true for $K=\bc$, it remains true if $K$ is replaced by $\bz$; then the same property allows us to conclude that it also holds for any algebraically closed field $K$ of characteristic zero.
\noindent So let us prove the claim (\ref{claimc}) for $K=\bc$.
As a consequence of Deligne's Correspondence Theorem (see \cite{BrSh}, Theorem 20.3.11) for all indices $i$ we have   
$$H_{J_{3,1}}^i(R)\simeq H^{i-1}(\bc^7\setminus V_{3,1}, \tilde R).$$ 
\noindent
Hence our claim can be restated equivalently as
$$H^3(\bc^7\setminus V_{3,1}, \tilde R)\ne0.$$
Therefore, in view of Lemma \ref{lemma8}, it suffices to show that 
\begin{equation}\label{2bis} H_{\rm DR}^{10}(\bc^7\setminus V_{3,1},\bc)\ne 0,\end{equation}
a statement that is the De Rham analogue to (\ref{2}) for $s=3$, $t=1$. For the sake of simplicity, we will omit the coefficient group $\bc$ in the rest of the proof. Let $W\subset K^7$ be the variety defined as in the proof of Theorem \ref{theorem5}, which in our present case is contained in $V_{3,1}$ and can be identified with the subvariety $V_{3,0}$ of $K^6$. By (\ref{6'}) we also have 
\begin{equation}\label{77}V_{3,1}\setminus W\simeq \bc^3\times(\bc\setminus\{0\}),\end{equation}
\noindent 
which is obviously non singular and pure-dimensional.
It is well known that
\begin{equation}\label{4bis}H^i(\bc^r)\simeq\left\{\begin{array}{cl} \bc&\mbox{if }i=0\\
0&\mbox{else, }
\end{array}\right.
\end{equation}
\noindent
and 
\begin{equation}\label{5bis}H^i(\bc^r\setminus\{0\})\simeq\left\{\begin{array}{cl} \bc&\mbox{if }i=0, 2r-1\\
0&\mbox{else. }
\end{array}\right.
\end{equation}
\noindent
Now, by  (\ref{77}) and the K\"unneth formula for singular cohomology (see \cite{W}, Theorem 3.6.1), 
$$H^i(V_{3,1}\setminus W)\simeq\displaystyle\bigoplus_{h+k=i}H^h(\bc^3)\otimes H^k(\bc\setminus\{0\}),$$
\noindent
so that, by (\ref{4bis}) and (\ref{5bis}), $H^i(V_{3,1}\setminus W)\ne0$ if and only if $i=0,1$. In particular
\begin{equation}\label{6bis} H^4(V_{3,1}\setminus W)= H^5(V_{3,1}\setminus W)=0.\end{equation}
Since, by (\ref{77}), the set $V_{3,1}\setminus W$ is a  closed non singular subvariety of $\bc^7\setminus W$ of codimension 3,
by \cite{Hart2}, Theorem 8.3, we have the following long exact sequence, which is the Gysin sequence for De Rham cohomology:
\begin{equation}\label{6ter}\cdots\rightarrow H_{\rm DR}^4(V_{3,1}\setminus W)\rightarrow H_{\rm DR}^{10}(\bc^7\setminus W)\rightarrow H_{\rm DR}^{10}(\bc^7\setminus V_{3,1})\rightarrow H_{\rm DR}^5(V_{3,1}\setminus W)\rightarrow 
 \cdots\end{equation}
\noindent
Now De Rham cohomology coincides with singular cohomology on non singular varieties, by virtue of Grothendieck's Comparison Theorem (see \cite{G}, Theorem 1', or \cite{Hart2}, Theorem, p. 147). Therefore, from (\ref{6bis}) it follows that the leftmost and the rightmost terms in (\ref{6ter}) vanish. Consequently,
\begin{equation}\label{7bis} H_{\rm DR}^{10}(\bc^7\setminus W) \simeq  H_{\rm DR}^{10}(\bc^7\setminus V_{3,1}) .
\end{equation}
\noindent
In view of  (\ref{7bis}), our claim (\ref{2bis}) will follow once we have proven that 
\begin{equation}\label{8bis}  H_{\rm DR}^{10}(\bc^7\setminus W)\ne0.\end{equation}
\noindent
 This is what we are going to show next. Recall from (\ref{ast}) and (\ref{decomposition}) that $\bc^7\setminus W=(\bc^6\setminus V_{3, 0})\cup Z$, where the union is disjoint and
\begin{equation}\label{decompositionbis}Z=\bc^6\times (\bc\setminus\{0\})\end{equation}
\noindent is a open subset of $\bc^7\setminus W$. We thus have the following Gysin sequence of De Rham cohomology:
\begin{equation}\label{9bis}\cdots\rightarrow H_{\rm DR}^9(Z)\rightarrow  H_{\rm DR}^8(\bc^6\setminus V_{3,0})\rightarrow H_{\rm DR}^{10}(\bc^7\setminus W)\rightarrow  H_{\rm DR}^{10}(Z)\rightarrow\cdots,\end{equation}
\noindent
where by  the K\"unneth formula for singular cohomology, (\ref{4bis}), (\ref{5bis}) and (\ref{decompositionbis}),
$$ H_{\rm DR}^{9}(Z)= H_{\rm DR}^{10}(Z)=0.$$
\noindent
 It follows that (\ref{9bis}) gives rise to an isomorphism: 
$$H_{\rm DR}^8(\bc^6\setminus V_{3,0})\simeq H_{\rm DR}^{10}(\bc^7\setminus W).$$
\noindent
But from \cite{BS}, Lemma 2, we know that $H^8(\bc^6\setminus V_{3,0})\ne0$,  so that $H_{\rm DR}^{10}(\bc^7\setminus W)\ne0$. This proves our claim (\ref{8bis}), which implies (\ref{2bis}) and shows claim (b) for $s=3$, $t=1$.  
\newline
Now suppose that $s>3$ and $t=1$. We have
$$J_{s,1}=I_s+(x_1z_1,\dots, x_sz_1).$$
\noindent
Hence
\begin{eqnarray}
J_{s,1}+(z_1)&=&I_s+(z_1),\label{first}\\
(J_{s,1})_{z_1}&=&(x_1,\dots, x_s)R_{z_1}.\label{second}
\end{eqnarray}
\noindent
We recall from \cite{BS}, Corollary, that 
\begin{equation}\label{a}{\rm cd}\,I_s=2s-3.\end{equation}
Note that the indeterminate $z_1$ does not occur in the minors generating $I_s$. Therefore, by virtue of Lemma \ref{plusone} we have  
${\rm cd}\,(I_s+(z_1))=\,{\rm cd}\,I_s+1.$  
In view of (\ref{first}) and (\ref{a}) it then follows that 
\begin{equation}\label{b}{\rm cd}\,(J_{s,1}+(z_1))=2s-2.\end{equation}
\noindent Moreover, since $x_1/1,\dots, x_s/1$  form a regular sequence in $R_{z_1}$,  they generate an ideal of cohomological dimension $s$ in  $R_{z_1}$. Thus, in view of (\ref{second}), we have
\begin{equation}\label{36'}{\rm cd}\,(J_{s,1})_{z_1}=s,\end{equation} 
\noindent where this cohomological dimension refers to the ring $R_{z_1}$. We have the following Brodmann sequence:
\begin{equation}\label{local}
\cdots\rightarrow H^{i-1}_{(J_{s,1})_{z_1}}(R_{z_1})\rightarrow H^i_{J_{s,1}+(z_1)}(R)\rightarrow H^i_{J_{s,1}}(R)\rightarrow H^i_{(J_{s,1})_{z_1}}(R_{z_1})\rightarrow\cdots,\end{equation}
\noindent
where we have used the fact that, due to the independence of base property (see \cite{BrSh}, Theorem 4.2.1, or \cite{Hu}, Proposition 2.11 (2)), $H^i_{J_{s,1}}(R_{z_1})\simeq H^i_{(J_{s,1})_{z_1}}(R_{z_1})$. Moreover, by (\ref {b}) and (\ref{36'}), 
$$H^i_{J_{s,1}+(z_1)}(R)=H^i_{(J_{s,1})_{z_1}}(R_{z_1})=0\qquad\mbox{for }i\geq 2s-1,$$
\noindent
because $s<2s-1$. Thus, in view of (\ref{local}), 
$$H^i_{J_{s,1}}(R)=0\qquad\mbox{for }i\geq2s-1.$$
\noindent
 We conclude that cd\,$J_{s,1}\leq 2s-2$. On the other hand, by (\ref{36'}), 
$$H^{2s-3}_{(J_{s,1})_{z_1}}(R_{z_1})=H^{2s-2}_{(J_{s,1})_{z_1}}(R)=0,$$
\noindent
because $s<2s-3$. Thus from (\ref{b}) and (\ref{local}) we deduce
$$0\neq H^{2s-2}_{J_{s,1}+(z_1)}(R)\simeq H^{2s-2}_{J_{s,1}}(R),$$
\noindent
which proves that cd\,$J_{s,1}\geq 2s-2$, whence we obtain cd\,$J_{s,1}=2s-2$, as claimed. \newline
Up to know we have proven claim (b) for all $s\geq 3$ and $t=1$, which settles the basis of our induction. Now we perform the induction step by assuming that $s\geq 3$, $t\geq 2$ and supposing that 
\begin{equation}\label{0}{\rm cd}\,J_{s,t-1}=2s+t-4.\end{equation} 
\noindent
We have 
\begin{eqnarray}
J_{s,t}+(z_t)&=&J_{s,t-1}+(z_t),\label{alpha}\\
(J_{s,t})_{z_t}&=&(x_1,\dots, x_s, y_1,\dots, y_{t-2})R_{z_t}\label{beta}.
\end{eqnarray}
Since the indeterminate $z_t$ does not occur in the generators of $J_{s,t-1}$ and the elements $x_1/1,\dots, x_s/1, y_1/1,\dots, y_{t-2}/1$ form a regular sequence in $R_{z_t}$, in view of Lemma \ref{plusone},  the relations (\ref{0}), (\ref{alpha}) and (\ref{beta})  allow us to deduce that
\begin{eqnarray}
{\rm cd}\,(J_{s,t}+(z_t))&=&2s+t-3\label{gamma}\\
{\rm cd}\,(J_{s,t})_{z_t}&=&s+t-2\label{delta},
\end{eqnarray}
where $s+t-2\leq 2s+t-3$. 
We have the following Brodmann sequence:
\begin{equation}\label{I}
\cdots\rightarrow H^{i-1}_{(J_{s,t})_{z_t}}(R_{z_t})\rightarrow H^i_{J_{s,t}+(z_t)}(R)\rightarrow H^i_{J_{s,t}}(R)\rightarrow H^i_{(J_{s,t})_{z_t}}(R_{z_t})\rightarrow\cdots.
\end{equation}
In view of (\ref{gamma}) and (\ref{delta}), in (\ref{I}) we have 
$$H^i_{J_{s,t}+(z_t)}(R)=H^i_{(J_{s,t})_{z_t}}(R_{z_t})=0\qquad\mbox{for }i\geq 2s+t-2,$$
\noindent 
which implies that cd\,$J_{s,t}\leq 2s+t-3.$ Moreover, from (\ref{delta}) we obtain 
$$H^{2s+t-4}_{(J_{s,t})_{z_t}}(R_{z_t})=H^{2s+t-3}_{(J_{s,t})_{z_t}}(R_{z_t})=0,$$
\noindent 
since $s>2$ implies that $2s+t-4>s+t-2$.
\noindent
Therefore, in view of (\ref{gamma}), in (\ref{I}) we have 
$$0\neq H^{2s+t-3}_{J_{s,t}+(z_t)}(R)\simeq H^{2s+t-3}_{J_{s,t}}(R),$$
\noindent
This yields cd\,$J_{s,t}=2s+t-3$, as claimed, and completes the proof.
\begin{remark}{\rm Theorem \ref{theorem5} and Theorem \ref{theorem9} (a) show that the inequality (\ref{cdara}) is strict for $J_{s,t}$ if $s\geq 3$ and char\,$K>0$. In fact, in this case we have 
$${\rm cd}\,J_{s,t}=s+t-1<s+t+s-3=2s+t-3\leq\,{\rm ara}\,J_{s,t}.$$
\noindent
According to Theorem \ref{theorem3}, however, equality  always holds for $s=2$; in turn, Theorem \ref{theorem9} (b) and Example \ref{Example1} show that  equality also holds for $s=3,4,5$ and $t=1$ provided that char\,$K=0$. The question in the remaining cases is open. In any case, Theorem \ref{theorem7} tells us that in characteristic zero the cohomological dimension and the arithmetical rank are close to each other, since
$${\rm ara}\,J_{s,t}\leq\,{\rm cd}\,J_{s,t}+1.$$ 
}
\end{remark}

\end{document}